\documentclass[11pt]{amsart}

\usepackage{amsmath, amssymb}
\usepackage{amsfonts}
\usepackage{amsthm}
\usepackage{amsopn}
\usepackage{amscd}

\usepackage [all]{xy}

\newcommand{\sC}{\mathcal{C}}
\newcommand{\cinf}{\sC^\infty}

\newcommand{\bR}{\mathbf{R}}

\newcommand{\bZ}{\mathbf{Z}}
\newcommand{\A}{\mathbb{A}}

\newcommand{\ruu}{\bR^{1|1}}
\newcommand{\rou}{\bR^{0|1}}

\newcommand{\tensor}{\otimes}
\newcommand{\comp}{\circ}

\newcommand{\tic}{\tilde{c}}
\newcommand{\ti}{\tilde}

\newcommand{\vd}{\partial_\th+\th\partial_t}
\newcommand{\sman}{\mathbf{SM}}
\newcommand{\str}{\textup{str}}

\theoremstyle{plain}
\newtheorem{thm}{Theorem}

\theoremstyle{definition}

\theoremstyle{remark}

\let\a\alpha

\let\d\delta

\let\th\theta

\let\o\omega
\let\G\Gamma

\let\O\Omega

\let\na\nabla
\let\ra\rightarrow
\let\lra\longrightarrow
\let\ti\tilde

\begin{document}
\title{A geometric view of the Chern character}
\author{Florin Dumitrescu}
\date{\today}
\maketitle

\begin{abstract} In this note we show that the Chern character form of a superconnection is obtained via the parallel transport of the superconnection along superpaths, by restriction to the universal superpoint path.
\end{abstract}

\vspace{.2in}
Consider a $\bZ/2$-graded vector bundle $E$ over a manifold $M$ and let $\A$ be a superconnection on $E$, i.e. an odd first-order differential operator on $E$ that satisfies the graded Leibniz rule:
\[ \A:\O^{*}(M, E)\to\O^{*}(M, E) \ \ \ \ \  \A(\o s)=d\o\tensor s+(-1)^{deg\ \o}\o\tensor\A(s), \]
for $\o\in\O^{*}(M)$ and $s\in \G(M,E)$. This notion generalizes the concept of a connection on a vector bundle in the sense that $\A$ can be split as a connection $\na$ on $E$ (grading-preserving) and a linear part $A\in \O^{*}(M, End\ E)^{odd}$. Superconnections were introduced by Quillen in \cite{Q} in order to obtain a local representation of the Chern character that would be better suited for a local form of the family index theorem, see \cite{AS4}. This problem was eventually solved by Bismut \cite{Bi} making use of the so-called Bismut superconnection in the infinite-dimensional setting (which we do not address here).

This article builds up on an idea of Fei Han \cite{Han} to obtain the Chern character form of a connection from the parallel transport along superpaths in the basespace $M$ determined by the connection, as defined in \cite{D1}. Moreover, Han shows that by considering parallel transport along superloops in $M$, one obtains the Bismut-Chern character form of \cite{Bi1}, an equivariant form on $LM$ the loopspace of $M$ which satisfies a differential equation (see Definition 6.4 of \cite{GJP}) reproduced by parallel transport along superloops. When restricted to $M$ (i.e. constant loops of $M$) this form gives the ordinary Chern character. Fei Han's work on Bismut-Chern character was explained to me by Stephan Stolz.

Given a superconnection $\A$ on a $\bZ/2$-bundle over a manifold $M$, we defined in \cite{D1} a notion of parallel transport along (families of) superpaths $c:\ruu\times S \ra M$ (parametrized by arbitrary supermanifolds $S$) which is compatible with glueing of superpaths, and is invariant under conformal reparametrizations of superpaths in the inverse adiabatic limit. This is done as follows. First, we write $\A=\na+A$, with $\na=\A_{1}$ the connection part of the superconnection $\A$ (which shifts the grading by 1 in $\O^{*}(M, E)$) and $A\in\O^*(M, End\ E)^{odd}$ the linear part of the superconnection. For an arbitrary superpath $c$ in $M$ consider the diagram
\[ \xymatrix @R=1.5pc @C=3pc { E \ar[dd] & & c^*E \ar[ll] \ar[dl] \ar[dd] \\
& \pi^*E \ar[ul] \ar[dd] & \\
M & & \ruu\times S \ar'[l][ll]_{\ \ \ \ c} \ar@{-->}[dl]^-{\tilde{c}} \\
& \Pi TM \ar[ul]^\pi &  } \]
with $\tic$ a canonical lift of the path $c$ to $\Pi TM$, the ``odd tangent bundle" of $M$. Recall that the superpath $c$ can be viewed as a map $c:\bR\times S\to \Pi TM$, via the identification 
\[ \sman(\ruu\times S, M)\cong\ \sman(\bR\times S, \Pi TM). \]Then $\tic$ is defined as the composition 
\[\ruu\times S= \bR\times \rou\times S \lra \rou\times \Pi TM \lra \Pi TM, \]
of the map $c$ followed by the action map $T$ of $\rou$ on $\Pi TM$, which, infinitesimally is given by the odd derivation (vector field) $d$ of functions on $\Pi TM$ (which are differential forms on $M$).

Then \emph{parallel transport along $c$} is defined by {\it parallel} sections $\psi\in\G(c^*E)$ along $c$ which are solutions to the following differential equation 
\[ (c^*\na)_{D}\psi- (\tilde{c}^*A)\psi= 0.  \]
Here $D=\vd$ denotes the standard (right invariant) vector field on $\ruu$, which generates the (super) Lie algebra of the super Lie group $\ruu$ and whose square is
\[ D^{2}=\frac{1}{2}[D,D]=\partial_{t}, \]
the standard time translation vector on $\bR$.
See the standard reference \cite{DM} or Section 2 of \cite{D1} for a brief introduction to supermanifolds.

Our main result here is to show that this ``$1|1$-parallel transport'' obtained from a superconnection, when restricted to the ``$0|1$-parallel transport'', it reproduces the Chern character form of the superconnection.

\begin{thm} Let $\A$ be a superconnection on a $\bZ/2$-bundle $E$ over a manifold $M$. The $1|1$-parallel transport along the superpath
 given by the composition
\[ \ruu\times \Pi TM\lra \rou\times \Pi TM\lra M, \]
where the first map is given by the projection $\ruu\to\rou$ and the second map is the ``superpoint evaluation map'' of $M$, gives rise to the Chern character form of the superconnection $ch(\A)=\str(\exp(-\A^{2}))$.
\end{thm}

\noindent\emph{Note.} In the above, the supertrace $\str$ is the extension of the ordinary supertrace defined on $\bZ/2$-graded endomorphisms:
\[ \str:\O^{*}(M,E)\lra \O^{*}(M):\ \ \ \o\tensor A\longmapsto \o\ \str A. \]

\begin{proof} Let us begin by remarking that this is a local problem, so it can be reduced to the case of a trivial bundle $E=\underline{\bR^{p}\oplus\bR^{q}}$ over $M$. \\

\noindent We consider first the case when the connection part $\na$ of the superconnection $\A$ is {\it flat}, as the calculation becomes more transparent. In this case, the connection can be taken to be the trivial connection $d$ (there is a trivialization of the bundle $E$ in which the connection is given by the trivial one). %On the other side, the parallel transport does not depend on the trivialization. 
Note that the superpath $c$ in $M$ given by the composition $ev\comp p$ lifts to a superpath $\tic$ in $\Pi TM$ given by the composition $T\comp p$ where $T: \rou\times\Pi TM\ra \Pi TM$ denotes the left action of $\rou$ on $\Pi TM$ which on functions is given by
\[ \O^{*}(M)\to \O^{*}(M)[\th] \]
\[ f\mapsto f+(df)\th, \ \ \text{ for } f\in\O^{0}(M) \]
\[ \a\mapsto\a+(-1)^{deg\ \a} (d\a)\th, \ \ \text{ for } \a\in\O^{*}(M).\]
The second relation is obtained from the first one, by taking into account that the pushforward of the odd vector field $d$ along the action map $T$ is again $d$. (The exterior derivative on forms on $M$ is interpreted as an odd vector field or derivation on $\Pi TM$ as $\cinf(\Pi TM)=\O^{*}(M)$. This vector field squares to zero, giving rise to an $\rou$ action- see \cite{D1}, Section 2.6-  on $\Pi TM$, which is the map $T$.) Let us represent the relevant maps in the diagram
\[ \xymatrix{ M & & \rou\times \Pi TM \ar[ll]_{ev} \ar[dl]_T & & \ruu\times \Pi TM. \ar[ll]_{p} \ar@{-->}[dlll]^{\ti{c}} \ar@/_2pc/[llll]_c \\
& \Pi TM \ar[ul]^\pi & & & } \]
The pullback connection of the trivial connection $d$ along $c=(ev)p$ is still the trivial connection $d$. The parallel transport equation along $c$ is given by
\[ (c^{*}d)_{D}(\psi_{0}+\th\psi_{1})-(\tic^{*}A)(\psi_{0}+\th\psi_{1})=0. \]
As $(c^{*}d)_{D}=D$ and $\tic^{*}A= A-\th dA$, the equation becomes
\[ (\vd)(\psi_{0}+\th\psi_{1})-(A-\th dA)(\psi_{0}+\th\psi_{1})= 0 \]
which is equivalent to
\[\psi_{1}+\th(\partial_{t}\psi_{0})-A\psi_{0}+\th A\psi_{1}+\th (dA)\psi_{0}=0. \]
This gives rise to the system
$$\left\{
\begin{array}{l}
\psi_1-A\psi_0=0\\
\frac{\partial\psi_0}{\partial t}+A\psi_1+(dA)\psi_0=0.
\end{array} \right.$$
Combining the last two relations we get
\[ \frac{\partial\psi_0}{\partial t}+(dA+A\wedge A)\psi_0=0. \]
The solution at $t=0$ and $t=1$ defines an element in $\O^{*}(M, End E)$ given by $\exp(-\A^{2})$, as $\A^{2}=dA+A\wedge A$. Taking the supetrace of the endomorphism valued form gives the Chern character form of the superconnection. \\
%The case of a nontrivial connection is similar.

\noindent Let us consider now the general case of an arbitrary superconnection on a (trivial) $\bZ/2$-graded vector bundle. Namely, we have $\A=\na+A$ and we write the connection $\na=d+\o$, with $\o\in\O^{1}(M, End^{0}\ E)$. Denote by $\d$ the exterior derivative on $\Pi TM$, $\rou\times\Pi TM$ and $\ruu\times \Pi TM$ (as we hope no confusion will arise). The pullback of $\na$ via $\pi$ is 
\[ \pi^{*}\na= \d+\pi^{*}\o. \]
The pullback of $\na$ via $c=(ev)p=\pi Tp$ is
\[ c^{*}\na=\d+p^{*}T^{*}\pi^{*}\o. \]
If $\o=fdg$, with $f,g$ functions on $M$, then $\pi^{*}\o=f\d g$. Further,
\begin{eqnarray*}
T^*(f\d g) & = & T^*(f)\d T^*(g) \\
& = & (f+df\th)\d(g+dg\th)\\
& = & (f+df\th)(\d g+\d(dg)\th-dg\d\th) \\
& = & f\d g+ f\d(dg)\th -fdg\d\th- df\d g\th + dfdg\th \d\th\\
& = & f\d g+(f\d(dg)-df\d g)\th+(-fdg+ dfdg\th)\d\th.
\end{eqnarray*}
Therefore, 
\[ <c^{*}(fdg)=p^{*}T^{*}\pi^{*}(fdg), \vd>=-fdg+ dfdg\th \]
as $g$ and $dg$ are independent of $t$ and $\th$. We obtain the formula
\[ <c^{*}(\o),\vd>=-\o+(d\o)\th. \]
The parallel transport equation along $c$ is given by
\[ (c^{*}\na)_{D}(\psi_{0}+\th\psi_{1})-(\tic^{*}A)(\psi_{0}+\th\psi_{1})=0. \]
As $$<(c^{*}\na), D>=<\d+c^{*}(\o),D>=D-\o+(d\o)\th$$ and $\tic^{*}A= A-\th dA$, the equation can be written
\[ (\vd)(\psi_{0}+\th\psi_{1})-(\o-\th d\o)(\psi_{0}+\th\psi_{1})-(A-\th dA)(\psi_{0}+\th\psi_{1})= 0 \]
or
\[ (\vd)(\psi_{0}+\th\psi_{1})-(A'-\th dA')(\psi_{0}+\th\psi_{1})= 0 \]
where $A'=A+\o\in\O^{*}(M, End \ E)^{odd}$. As before, this is equivalent to
\[\psi_{1}+\th(\partial_{t}\psi_{0})-A'\psi_{0}+\th A'\psi_{1}+\th (dA')\psi_{0}=0. \]
which produces the system
$$\left\{
\begin{array}{l}
\psi_1-A'\psi_0=0\\
\frac{\partial\psi_0}{\partial t}+A'\psi_1+(dA')\psi_0=0.
\end{array} \right.$$
Combining these relations we get
\[ \frac{\partial\psi_0}{\partial t}+(dA'+A'\wedge A')\psi_0=0. \]
The solution at $t=0$ and $t=1$ defines an even element in $\O^{*}(M, End\ E)$ given by $\exp(-\A^{2})$, as $\A=\na+A= d+\o+A= d+A'$ and $\A^{2}=dA'+A'\wedge A'$. The supertrace gives us the Chern character form of the superconnection.

\end{proof}

\noindent\emph{Acknowledgements.}  This paper was written while the author was on a postdoctoral position at University of Hamburg. I would like to thank Tilmann Wurzbacher, whose invitation to speak at the seminar in Bochum gave me the opportunity to review some of the material.

\bibliographystyle{plain}
\bibliography{bibliografie}

\bigskip
\raggedright Universit\"at Hamburg\\  Bundesstra\ss{}e 55 \\
D-20146 Hamburg\\ Email: {\tt florinndo@gmail.com}

\end{document}